\documentclass[12pt,a4paper,english,intoc,bibliography=totoc,index=totoc,BCOR10mm,captions=tableheading,titlepage,fleqn]{scrbook}
\usepackage{lmodern}

\usepackage[T1]{fontenc}
\usepackage[latin9]{inputenc}
\usepackage{babel}
\usepackage{amsmath}
\usepackage{amssymb}
\usepackage[unicode=true,
 bookmarks=true,bookmarksnumbered=true,bookmarksopen=true,bookmarksopenlevel=2,
 breaklinks=false,pdfborder={0 0 0},backref=false,colorlinks=false]
 {hyperref}
\hypersetup{pdftitle={Abstract},
 pdfauthor={Uwe Stöhr},
 pdfpagelayout=OneColumn, pdfnewwindow=true, pdfstartview=XYZ, plainpages=false}
\usepackage{breakurl}

\makeatletter

\usepackage{calc}

\makeatother

\begin{document}
\noindent \begin{center}
\textbf{\Large Compact Sets of Baire Class One Functions and Maximal
Almost Disjoint Families}
\par\end{center}{\Large \par}

\noindent \begin{center}
{\large Haim Horowitz and Stevo Todorcevic}
\par\end{center}{\large \par}

\noindent \begin{center}
\textbf{\small Abstract}
\par\end{center}{\small \par}

\noindent \begin{center}
{\small We provide a proof that analytic almost disjoint families
of infinite sets of integers cannot be maximal using a result of Bougain
about compact sets of Baire class one functions. Inspired by this
and related ideas, we then provide a new proof of that there are no
maximal almost disjoint families in Solovay's model. We then use the
ideas behind this proof to provide an extension of a dichotomy result
by Rosenthal and by Bourgain, Fremlin and Talagrand to general pointwise
bounded functions in Solovay's model. We then show that the same conclusions
can be drawn about the model obtained when we add a generic selective
ultrafilter to the Solovay model.}
\par\end{center}{\small \par}

\textbf{\large Introduction}{\large \par}

The initial motivation for this paper is the study of the definability
of maximal almost disjoint families of infinite sets of integers but
we soon realsed that this could be seen as part of the more general
study of pointwise convergence of sequences of continuous functions
on separable metric spaces. This first line of study goes back to
the classical result of Mathias ({[}Ma{]}) showing that there are
no analytic mad families while the second goes back to the paper of
Bourgain ({[}Bo{]}) about pointwise convergence of sequences of continuous
functions on Polish spaces. Recently two more new proofs of Mathias'
results have been discovered, one by Toernquist and the other by Bakke
Haga, Schrittesser and Toernquist (see ({[}To{]}) and ({[}BST{]}).
In the first part of the paper, we shall provide another proof of
Marthias' result using a result by Bourgain ({[}Bo{]}) on Baire class
one functions. Motivated by that proof, we shall then provide a new
proof of the nonexistence of mad families in Solovay's model. Recall
that the study of the nonexistence of mad families in choiceless models
was established by Mathias in ({[}Ma{]}), where he proved that there
are no mad families in the Solovay model obtained from a Mahlo cardinal.
The upper bound on the consistency strength was later reduced to an
inaccessible cardinal by Toernquist ({[}To{]}) and to $ZFC$ by the
first author and Shelah ({[}HS{]}).

Our proof of the nonexistence of mad families in Solovay's model will
actually provide a more general dichotomy result involving general
selective coideals (see Theorem C in the last section), which will
have several interesting corollaries, such as an extension of a result
by Godefroy ({[}Go{]}) from analytic sets to general sets of reals
in Solovay's model, the fact that $\mathfrak s=\mathfrak c$ in Solovay's
model and an extension of a dichotomy result originally due to Bourgain,
Rosenthal and Bourgain-Fremlin-Talagrand:

\textbf{Theorem ({[}Bo{]}, {[}Ro{]}, {[}BFT{]}): }Let $(f_n : n<\omega)$
be a sequence of pointwise bounded continuous functions from a Polish
space $X$ to $\mathbb R$. Suppose that a Baire class one function
$f$ is in the pointwise closure of the sequence $(f_n : n<\omega).$
Then we have one of the following two alternatives: 

a. $(f_n : n<\omega)$ contains a subsequence whose closure is homeomorphic
to $\beta \omega$.

b. $(f_n : n<\omega)$ contains a subsequence converging pointwise
to $f$.

We shall finish the paper by showing that all these results about
the Solovay model $L(\mathbb{R})[\mathcal{U}]$ where $\mathcal{U}$
is a (generic) selective ultrafilter on $\omega.$ So, in particular,
we show that there are no mad families in $L(\mathbb{R})[\mathcal{U}]$
which is a considerably stronger from the clam that there are no such
families in the smaller model $L(\mathbb{R}).$

\textbf{\large Baire class one functions and analytic almost disjoint
families}{\large \par}

The main result of this section will be a new proof of the nonexistence
of analytic mad families where our main tool will be a result of Bourgain
on Baire class one functions {[}Bo{]}.

\textbf{Theorem A: }There are no analytic mad families.

\textbf{Proof: }Assume towards contradiction that there is an analytic
mad family, and fix such a family $\mathcal A$.

\textbf{Definition A.1: }Let $X_{\mathcal{A}}:=\{ \bold 0\} \cup \{ f_n : n<\omega\} \cup \{ \delta_x : x\in \mathcal A\}$
where:

a. $\bold 0 : \mathcal A \rightarrow \mathbb R$ is contantly $0$.

b. $f_n : \mathcal A \rightarrow \{0, 1\}$ is defined by $f_n(x)=1$
iff $n\in x$.

c. $\delta_x : \mathcal A \rightarrow \{0,1\}$ is defined by $\delta_x(y)=1$
iff $y=x$.

\textbf{Observation A.2: $X_{\mathcal A}$ }is sequentially compact.

\textbf{Proof: }Let $\{ g_n : n<\omega \} \subseteq X_{\mathcal A}$,
we may assume wlog that, for every $k<\omega$, $g_k=f_{n_k}$ for
some $n_k<\omega$. Let $x=\{ n_k : k<\omega \}$, we may assume wlog
that $x$ is infinite and the $n_k$ are pairwise distinct. By the
madness of $\mathcal A$, there is some $y\in \mathcal A$ such that
$|x\cap y|=\aleph_0$. Denote $x\cap y$ by $I$. For every $n_k \in I$,
$f_{n_k}(y)=1$. For every $y' \in \mathcal A$, if $y' \neq y$ then
$|y' \cap y|<\aleph_0$, hence for every $k$ large enough, if $n_k \in I$
then $n_k \notin y'$, hence $f_{n_k}(y')=0$. It follows that $\{f_{n_k} : n_k \in I\}$
converges pointwise to $\delta_y$, therefore, $X_{\mathcal A}$ is
sequentially compact. $\square$

\textbf{Observation A.3: }$\{ \bold 0\}$ is in the pointwise closure
of $\{f_n : n<\omega\}$. $\square$

\textbf{Observation A.4: $\{f_n : n<\omega\}$ }has no subsequence
that pointwise converges to $\bold 0$.

\textbf{Proof: }By a similar argument as before, i.e. suppose that
some subsequence $\{ f_{n_k} : k<\omega \}$ converges pointwise to
$\bold 0$ and let $x=\{ n_k : k<\omega\}$. By madness, there is
$y\in \mathcal A$ such that $|x\cap y|=\aleph_0$, and as before,
$\{f_{n_k} : n_k \in I\}$ converges pointwise to $\delta_y$, contradicting
our assumption. $\square$

We now use the following result by Bourgain:

\textbf{Theorem A.5 {[}Bo{]}: }Let $Y$ be a Polish space and let
$H$ be a subset of $B_1(Y)$. If $H$ is relatively compact in $B_1(Y)$,
then any limit point of a sequence $\{g_n : n<\omega\}$ in $H$ is
the pointwise limit of a subsequence of $\{g_n : n<\omega\}$.

\textbf{Proof: }This is essentially the content of the proof of Theorem
12 in {[}Bo{]}. $\square$

\textbf{Observation A.6: }As $\mathcal A$ is analytic, there is a
continuous surjection $g: \mathbb R \rightarrow \mathcal A$, which
naturally induces an embedding $G: X_{\mathcal A} \rightarrow B_1(\mathbb R)$.
Furthermore, observations 2-4 hold for the images of $X_{\mathcal A}$,
$\{ f_n : n<\omega\}$ and $\bold 0$ under $G$, and we shall identify
these objects with their images under $G$. $\square$

\textbf{Observation A.7: }$X_{\mathcal A}$ is relatively compact
in $B_1(\mathbb R)$.

\textbf{Proof: }By Theorem 4 in {[}Bo{]} and the fact that $X_{\mathcal A}$
is sequentially compact. $\square$

Finally, by Observtion 7, Observation 3 and Theorem 5, there is a
subsequence of $\{f_n : n<\omega\}$ that pointwise converges to $\bold 0$.
By Observation 4, we get a contradiction. This completes the proof
of the theorem. $\square$

$\\$

\textbf{\large Mad families and selective coideals in Solovay's model}{\large \par}

Inspired by the ideas from the previous section, the main result of
this section will be a new proof of the nonexistence of mad families
in Solovay's model. We work with the same characteristic functions
$(f_n : n<\omega)$ from the previous section. We shall enumerate
the subsets $X_n \subseteq [\omega]^{\omega}$ consisting of generic
branches through trees in $V$. Given a candidate $\mathcal A$ for
a mad family in Solovay's model and the derived coideal $\mathcal H$
of sets that are not almost covered by finite unions from $\mathcal A$,
we derive a subset $Y\subseteq \mathcal A$ and $M\in \mathcal H$
such that there is no $M' \subseteq M$ from $\mathcal H$ such that,
for some $l<\omega$, $X_l \cap Y \neq \emptyset$ and $f_n \restriction X_l$
is contsant for every $n\in M'$. A similar argument to that establishing
the perfect set property in Solovay's model will show that there is
some infinite $X=\{k_n : n<\omega\} \subseteq M$ such that $\mathcal A \restriction X=\{Y\cap X : Y\in \mathcal A\}=\mathcal{P}(X)$,
this will contradict the almost disjointness of $\mathcal A$.

\textbf{Theorem B: }There are no mad families in Solovay's model.

\textbf{Proof: }Assume towards contradiction that the theorem fails. 

Let $\kappa$ be an inaccessible cardinal and let $G\subseteq Col(\omega,\kappa)$
be generic over $V$. Let $\mathcal A \in V[G]$ be a mad family definable
from some $a\in Ord^{\omega}$ using the formula $\phi$ and define
$(f_n : n<\omega)$ and $\bold 0$ as before. Observations A.2 and
A.4 only use the maximality of $\mathcal A$, hence they hold in this
context as well. 

In $V[G]$, there exists an enumeration $(T_n : n<\omega)$ of all
subtrees of $2^{<\omega}$ from $V$. For $n<\omega$, let $[T_n]^{V[G]}=\{ x\in 2^{\omega} \cap V[G] : x$
is a branch of $T_n\}$ and let $X_n \subseteq \mathcal{P}(\omega)$
be the set $\{x^{-1}(1) : x\in [T_n]^{V[G]} \}$.

Let $\mathcal H$ be the coideal of all sets $X\subseteq \omega$
such that $\bold 0$ is an accumluation point of $\{f_n : n\in X\}$.
Note that $\mathcal H$ is simply the coideal of sets that are not
covered (modulo a finite set) by a finite union of elements of $\mathcal A$.
As $\mathcal A$ is mad, every decreasing $\omega$-sequence of elements
of $\mathcal H$ has a pseudointersection in $\mathcal H$. We shall
construct a sequence $((M_n,Y_n) : n<\omega)$ by induction on $n<\omega$
as follows:

\textbf{Case I, $n=0$: }Let $Y_0=\mathcal A$ and $M_0$ be any member
of $\mathcal H$.

\textbf{Case II, $n=k+1$: }If there is a $l<\omega$ such that $X_l \cap Y_k \neq \emptyset$
and there is $M_* \subseteq M_k$ such that $M_* \in \mathcal H$
and $f_m$ is contant on $X_l$ for every $m\in M_*$, then we let
$M_n=M_*$ and $Y_n=Y_k \setminus X_l$ where $l$ is the least natural
number with this property. If there are no such $l$ and $M'$, we
stop the induction.

\textbf{Definition: }We shall define $(M,Y)$ as follows:

a. If there is $n<\omega$ for which we can't proceed to the $n+1$th
stage, we let $(M,Y)=(M_n,Y_n)$.

b. If there is no such $n<\omega$, i.e. if we carried the induction
successfully, we let $M \in \mathcal H$ be a pseudo intersection
of $(M_n : n<\omega)$ and $Y=\underset{n<\omega}{\cap}Y_n$. 

\textbf{Claim: }$Y$ is uncountable in $V[G]$.

\textbf{Proof: }Suppose towards contradiction that $Y=\{x_n : n< \omega \}$
is countable. Therefore, wlog $\mathcal A \subseteq \underset{n<\omega}{\cup}X_n \cup \{x_n : n<\omega\}$.
As $M\subseteq^* M_1$, there is $n_1<\omega$ such that $f_m$ is
constant on $X_0$ for every $m\in M \setminus n_1$. Therefore, there
is $M_1' \subseteq M$ and $i_1 \in \{0,1\}$ such that $M_1' \in \mathcal H$
and $f_m(x)=i_1$ for every $m\in M_1'$ and $x\in X_0$. Similarly,
$M_1' \subseteq^* M_2$, hence there are $M_2' \subseteq M_1'$ and
$i_1 \in \{0,1\}$ such that $M_2' \in \mathcal H$ and $f_m(x)=i_1$
for every $m\in M_2'$ and $x\in X_1$. We proceed in a similar way
by induction on $n<\omega$, obtaining sets $M_n' \in \mathcal H$
and numbers $i_n \in \{0,1\}$. Finally, let $M' \in \mathcal H$
such that $M' \subseteq^* M_n'$ for every $n < \omega$. Now let
$N_0 \subseteq M'$ be a set from $\mathcal H$ and let $j_0 \in \{0, 1\}$
such that $f_m(x_0)=j_0$ for every $m\in N_0$. Proceed by induction
to obtain a decreasing sequence $(N_n : n<\omega)$ and $(j_n : n<\omega)$
such that $N_n \in \mathcal H$ and $f_m(x_n)=j_n$ for every $m\in N_n$.
Finally, let $N \in \mathcal H$ be a set such that $N \subseteq^* N_n$
for every $n<\omega$. Now define $f: \mathcal A \rightarrow \{0,1\}$
in the following way: for every $x\in \mathcal A$, the value of $f_m(x)$
becomes constant for large enough $m\in N$, and we define $f(x)$
to be that value. Therefore, $(f_n : n\in N)$ pointwise converges
to $f$. As in the proof of Observation A.4, letting $y\in \mathcal A$
such that $|y\cap N|=\aleph_0$, the sequence $(f_n : n\in N\cap y)$
converges pointwise to $\delta_y$. It follows that $(f_n : n\in N)$
converges pointwise to $\delta_y$. Therefore, there is a set in $\mathcal H$
that doesn't have $\bold 0$ as an accumulation point, contradicting
the definition of $\mathcal H$. This completes the proof of the claim
as it follows that $Y$ is uncountable. $\square$

By the construction of $Y$, it's definable from some real in $V[G]$
and wlog we may assume that $Y$ and $\mathcal A$ are both definable
from $a$. Let $\alpha<\kappa$ such that $a\in V[G\cap Col(\omega,\alpha)]$.
As $Y$ is uncountable in $V[G]$, there exists some $y\in Y \setminus V[G\cap Col(\omega,\alpha)]$.
There exist $\beta \in (\alpha, \kappa)$, a $Col(\omega,\beta )$-name
$\underset{\sim}{y}$ and a condition $p_0 \in Col(\omega,\beta) \cap G$
such that $y=\underset{\sim}{y}[G\cap Col(\omega,\beta )]$ and (in
$V[G\cap Col(\omega,\alpha)]$) $p_0 \Vdash_{Col(\omega,\beta )} " \underset{\sim}{y} \in \underset{\sim}{Y}"$
(i.e. $\underset{\sim}{y}$ satisfies the formula that defines $Y$). 

Let $\mathbb P=Col(\omega,\beta)$ and $(D_n : n<\omega)$ be an enumeration
in $V[G]$ of the dense subsets of $\mathbb P$ from $V[G\cap Col(\omega,\alpha)]$.
We shall now construct a perfect tree $(p_t : t\in 2^{<\omega})$
of conditions of $\mathbb P$ and a Cantor scheme $(U_t : t\in 2^{<\omega})$
by induction on $|t|$ as follows:

\textbf{Case I, $|t|=0$: }Choose some $p_{()}\in D_0$ above $p_0$
and let $U_{()}$ be the set all reals $x$ such that $x(m)=\underset{\sim}{y}(m)$
whenever $p_{()}$ decides $\underset{\sim}{y}(m)$. 

\textbf{Case II, $|t|=n+1$: }Let $(s_i : i\in 2^n)$ list $2^n$
and suppose that $(p_{s_i} : i\in 2^n)$ and $(U_{s_i} : i\in 2^n)$
were chosen such that $p_{s_i} \in D_n$ and $U_{s_i}$ is the set
of all reals $x$ such that $x(m)=\underset{\sim}{y}(m)$ whenever
$p_{s_i}$ decides $\underset{\sim}{y}(m)$. Now suppose towards contradiction
that for every $m\in M$ there is some $i\in 2^n$ such that $p_{s_i}$
decides $\underset{\sim}{y}(m)$. Denote the set of $m\in M$ such
that $\underset{\sim}{y}(m)$ is decided by $p_{s_i}$ by $M_i$,
then by our assumption, $M=\underset{i\in 2^n}{\cup}M_i$. Therefore,
there is some $i\in 2^n$ such that $M_i \in \mathcal H$. Let $T \in V$
be the subtree of $2^{<\omega}$ determinned by $p_{s_i}$, i.e.,
given $t\in T$ of length $n$, if $p_{s_i} \Vdash "\underset{\sim}{y}(n)=j"$,
then $t \hat{} (j)$ is the only successor of $t$ in $T$. Otherwise,
both possible successors of $t$ are in $T$. In $V[G]$, there is
a set $H\subseteq \mathbb P$ such that $p_{s_i} \in H$ and $H$
is generic over $V[G\cap Col(\omega,\alpha)]$. Therefore, $\underset{\sim}{y}[H] \in \{ x^{-1}(1) : x\in [T]^{V[G]} \} \cap Y$
and for each $m\in M_i$, $f_m$ is constant on $\{x^{-1}(1) : x\in [T]^{V[G]} \}$,
which has the form $X_n$ for some $n<\omega$. We now obtain a contradiction
to the choice of $(M,Y)$, as we can find a subset $M' \subseteq M$
and a set $X_n$ such that $M' \in \mathcal H$, $X_n \cap Y \neq \emptyset$
and $f_m$ is constant on $X_n$ for every $m\in M'$. 

It follows that there is some $k \in M$ such that $\underset{\sim}{y}(k)$
is not decided by any of the conditions $(p_{s_i} : i\in 2^n)$, let
$k_n$ be the minimal such $k$. For every $s_i \in 2^n$, let $p_{s_i \hat{} 0}$
and $p_{s_i \hat{} 1}$ be two extensions of $p_{s_i}$ in $D_{n+1}$
such that $p_{s_i \hat {} l} \Vdash "\underset{\sim}{y}(k_n)=l"$
$(l=0,1)$. Finally, for $t\in 2^{n+1}$, let $U_t$ be the set of
all elements $x$ such that $x(m)=\underset{\sim}{y}(m)$ whenever
$p_t$ decides $\underset{\sim}{y}(m)$. If necessary, we may increase
$p_{s_i \hat {} l}$ inside $D_{n+1}$ to guarantee that $diam(U_{s_i \hat{} l}) \leq 2^{-(n+1)}$.

Now let $\{k_n : k<\omega\} \subseteq M$ be the sequence of $k_n$s
constructed during the induction. We may assume wlog that the sequence
is strictly increaing (in the proof of the existence of $k_n$, we
can replace $M$ by $M\setminus k_{n-1}$). For every $x\in 2^{\omega}$,
consider the set $G_x:=\{q \in \mathbb{P} :$ there exists some $t\leq x$
such that $q\leq p_t \}$. As $p_t \in D_{|t|}$ for every $t\in 2^{<\omega}$,
it follows that $G_x$ is $\mathbb P$-generic over $V[G\cap Col(\omega,\alpha)]$.
Therefore, $\underset{\sim}{y}[G_x] \in Y$ (as $p_0$ forces this)
and $\underset{\sim}{y}[G_x] \in \underset{n<\omega}{\cap}U_{x\restriction n}$.
Now let $A\neq B \subseteq \{k_n : k<\omega\}$ be infinite subsets
with an infinite intersection. Define $x_A \in 2^{\omega}$ by $x_A(n)=1$
iff $k_n \in A$ and define $x_B$ similarly. Let $y_A=\underset{\sim}{y}[G_{x_A}]$
and let $y_B=\underset{\sim}{y}[G_{x_B}]$, then $y_A,y_B \in \mathcal A$.
By the choice of the conditions $p_t$, $k_n \in y_A \cap y_B$ iff
$k_n \in A\cap B$, therefore, $|y_A \cap y_B|=\aleph_0$ and $y_A \neq y_B$,
contradicting the almost disjointness of $\mathcal A$. This completes
the proof of Theorem B. $\square$

$\\$

\textbf{\large Corollaries of Theorem B}{\large \par}

In this section we shall derive a few quick corollaries from the proof
of Theorem B. We first note that the above proof actually proves the
following, more general result:

\textbf{Theorem C: }Let $\mathcal H \subseteq \mathcal{P}(\omega)$
be a selective coideal and let $\mathcal A$ be a set of reals in
Solovay's model, then one of the following holds:

a. For every $M_0 \in \mathcal H$, there is an infinite set $X \subseteq M_0$
such that $\mathcal A \restriction X=\{Y\cap X : Y \in \mathcal A \}=\mathcal{P}(X)$.

b. There is $M\in \mathcal H$ such that $(f_n : n\in M)$ converges
pointwise to some function $f$ where $(f_n : n<\omega)$ are as before.
$\square$

Combining the above theorem with clause (b'') of the observation below
will provide an extension of an older result by Godefroy that was
previously established for analytic sets ({[}Go{]}).

\textbf{Observation: }Clause (b) in Theorem C is equivalent to:

b'. For every $x\in \mathcal A$, $M\subseteq^* x$ or $M\subseteq^* \omega \setminus x$.

b''. Every element in $\mathcal A \restriction M$ is either a finite
or a cofinite subset of $M$. $\square$

\textbf{Corollary: }Let $\mathcal A$ be a splitting family in Solovay's
mode, then there is an infinite $X\subseteq \omega$ such that $\mathcal A \restriction X=\mathcal{P}(X)$,
hence $|\mathcal{A}|=|\mathbb R|$ in Solovay's model. It folllows
that $\mathfrak s=\mathfrak c$ in Solovay's model. $\square$

In thee proof of Theorem B (which, as noted, is also a proof of Theorem
C), for each $x\in 2^{\omega}$, we let $A_x$ be the unique member
of $\underset{n<\omega}{\cap}U_{x\restriction n}$. Let $P=\{ A_x : x\in 2^{\omega} \}$,
obviously, $P$ is homeomorphic to $2^{\omega}$ and $f_{k_n}(A_x)=1$
iff $x(n)=1$. Now observe that the closure of $\{f_{k_n} : n<\omega\}$
in $2^P$ is homeomorphic to $\beta \omega$: $\beta \omega$ is the
closure in $2^{2^{\omega}}$ of $\{f_n : n<\omega\}$. By the above
remark, this is homeomorphic to the closure of $\{f_{k_n} \restriction P : n<\omega\}$
in $2^P$. Also note that the last argument is valid for general functions
in $\{0, 1\}$ (not just functions as in Definition A.1). The following
corolary now follows:

\textbf{Corollary: }Let $\mathcal A$ be a set of reals in Solovay's
model and let $\{f_n : n<\omega \}$ be a set of functions in Solovay's
model from $\mathcal A$ to $\{0, 1\}$. Given a selective coideal
$\mathcal H \subseteq \mathcal{P}(\omega)$, one of the following
holds:

a. For every $M_0 \in \mathcal H$ there is an infinite $X\subseteq M_0$
such that the closure of $\{f_n : n\in X\}$ has cardinality $>2^{\aleph_0}$.

b. There is some $M\in \mathcal H$ such that $(f_n : n\in M)$ is
pointwise convergent. $\square$

The above corollary will now imply the desired extension of the dichotomy
theorem by Rosenthal and by Bourgain-Fremlin-Talagrand to arbitrary
functions in Solovay's model from an arbitrary set of reals into $2^{\omega}$: 

\textbf{Theorem D: }Let $\mathcal A$ be a set of reals in Solovay's
model and let $(f_n : n<\omega)$ be a sequence of functions in Solovay's
model from $\mathcal A$ into  $2^\omega$, then one of the following
holds:

a. There is an infinite $X\subseteq \omega$ such that the closure
of $\{ f_n : n\in X \}$ has cardinality $>2^{\aleph_0}$.

b. $(f_n : n<\omega)$ has a converging subsequence.

\textbf{Proof: }
We shall use the previous corollary for the selective
coideal $\mathcal H_0=[\omega]^{\omega}$. For $n<\omega$, we define
the functions $(f_n^k : k<\omega)$ by $f_n^k(x)=f_n(x)(k)$. Consider
the sequence $\bar{f}_0=(f_n^0 : n<\omega)$, this sequence satisfies
the assumptions of the previous corollary. If clause (a) of the previous
corollary holds for $\bar{f}_0$, then we are done. Otherwise, there
is $M_0 \in \mathcal H_0$ as in clause (b). Now consider the sequence
$\bar{f}_1=(f_n^1 : n\in M_0)$ and the selective coideal $\mathcal{H}_1=[M_0]^{\omega}$.
If clause (a) of the previous corollary holds for $\bar{f}_1$, then
we are done. Otherwise, there is $M_1 \subseteq M_0$ as in clause
(b) for $\bar{f}_1$. We shall continue following this procedure similarly
for every $m<\omega$. If for some $m<\omega$ clause (a) of the previous
corollary holds, then we're done. Otherwise, we shall obtain a sequence
$(M_m : m<\omega)$ of infinite sets such that $M_{m+1} \subseteq M_m$
for every $m<\omega$. Now let $M \in [\omega]^{\omega}$ be a pseudointersection
of $\{M_m : m<\omega\}$, then obviously $(f_m : m \in M)$ is pointwise
convergent, which completes the proof. $\square$

Let $\{(p_i, q_i) : i<\omega\}$ be an enumeration of all ordered
pairs $(p, q)$ or rational numbers such that $p<q$. For $i<\omega$,
let $f_n^i(x)$ be defined as $0$ if $f_n(x) \leq p_i$, $1$ if
$q_i \leq f_n(x)$ and $2$ otherwise. For $i<\omega$, we define
the sequences $\bar{f_i}$ as in the previous proof and repeat the
same argument as before. It can be then shown that:

\textbf{Theorem E: }In Theorem D, we can further assume that the functions
$f_n$ are into $\mathbb R$. $\square$

\textbf{Theorem F (Solovay's model): }Let $\{f_n : n<\omega \}$ be
continuous functions from a set of reals $X$ to $\mathbb R$ and
let $f$ be in the closure of $\{ f_n : n<\omega\}$, then one of
the following hold:

a. $\mathcal{H}_f := \{N \subseteq \omega : f$ is in the closure
of $\{f_n : n\in N\} \}$ is a selective coideal.

b. There is a perfect set $P\subseteq X$ and a subsequence $\{f_{k_n} : n<\omega\} \subseteq \{f_n : n<\omega\}$
that behave as in the proof of Theorem B (i.e. the $f_{k_n}$ behave
like projections).

\textbf{Proof: }Let $W$ be the closure of $\{ f_n : n<\omega\}$.
By Lemma 2 in Section 12 of {[}Tod{]}, if $W$ is countably tight,
then $\mathcal{H}_f$ is selective, so assume that $W$ is not countably
tight. Therefore, there is $Z\subseteq$ and some $g$ in the closure
of $Z$ such that $g$ is not in the closure of any countable $A\subseteq Z$.
By Corollary 4 in Section 10 of {[}Tod{]}, we may assume wlog that
$g$ is the zero function and all functions in $Z$ are positive.
Given a countable $A\subseteq Z$ and $\epsilon>0$, let $X_{\epsilon}(A)=\{x\in X : \epsilon \leq f(x)$
for all $f\in A\}$. Suppose that for every $\epsilon>0$ there is
a countable $A_{\epsilon}\subseteq Z$ such that $X_{\epsilon}(A_{\epsilon})=\emptyset$,
then $g$ is in the closure of $\underset{n<\omega}{\cup}A_{\frac{1}{n}}$,
contradicting our assumption. Therefore, suppose that there is an
$\epsilon>0$ without the above property (which will be fixed until
the end of the proof), and we shall prove that clause (b) of the theorem
holds. For a countable $A\subseteq Z$, let $X(A)$ be $X_{\epsilon}(A)$.
By our assumption, $X(A) \neq \emptyset$ for all countable $A\subseteq Z$.
Note also that $A\subseteq B \rightarrow X(B) \subseteq X(A)$. We
shall try to construct an increasing sequence $(A_{\alpha} : \alpha<\omega_1)$
by induction on $\alpha<\omega_1$ as follows: $A_0$ will be any
countable subset of $Z$. If $\delta<\omega_1$ is a limit ordinal,
$A_{\delta}=\underset{\beta<\delta}{\cup}A_{\beta}$. At a successor
stage $\alpha+1$, we consider the sets $(X_n : n<\omega)$ from the
beginning of the proof of Theorem B. If there is a countable $A\subseteq Z$
and $n<\omega$ such that $A_{\alpha} \subseteq A$, $X(A_{\alpha}) \cap X_n \neq \emptyset$
and $X(A) \cap X_n=\emptyset$, we let $A_{\alpha+1}=A$. The process
must stop at a countable succesor ordinal $\alpha+1$, and we let
$B=A_{\alpha+1}$. We shall now consider $X(B)$. Note that $X(B)$
is uncountable: Suppose not, then $X(B)=\{ x_n : n<\omega \}$. For
each $n<\omega$, choose some $g_n \in Z$ such that $g_n(x_n)<\epsilon$.
Now note that $X(B\cup \{g_n : n<\omega\})=\emptyset$, contradicting
the assumption on $\epsilon$. It follows that $X(B)$ is uncountable.
Let $\underset{\sim}{y}$ be a name for an element in $X(B)$ as in
the proof of Theorem B. We let $\mathbb P$ and $(D_n : n<\omega)$
be as in the proof of Theorem B, and we shall construct by induction
on $|t|$ a tree of conditions $(p_t : t\in 2^{<\omega})$ and a scheme
$(U_t : t\in 2^{<\omega})$ as there. At stage $n+1$, we let $(p_{s_i} : i \in 2^n)$
and $(U_{s_i} : i\in 2^n)$ be as in the proof of Theorem B. For each
$i\in 2^n$, pick some $x_i \in U_{s_i} \cap X(B)$. There is some
$g_i \in Z$ such that $g_i(x_i)<\frac{\epsilon}{4}$. By the choice
of $B$, there is some $y_i \in U_{s_i} \cap X(B)$ such that $\epsilon \leq g_i(y_i)$.
As $g_i$ is in the closure of $\{f_n : n<\omega\}$, there is a large
enough $k_n$ such that $f_{k_n}(x_i)<\frac{\epsilon}{3}$ and $\frac{2}{3 \epsilon}< f_{k_n}(y_i)$.
By continuity, there is a large enough $j(n)$ such that for all $i\in 2^n$,
$f_{k_n}(x)<\frac{\epsilon}{3}$ if $x\restriction j(n)=x_i \restriction j(n)$
and $\frac{2}{3 \epsilon}<f_{k_n}(y)$ if $y\restriction j(n)=y_i \restriction j(n)$.
For each $i\in 2^n$, let $p_{s_i \hat{} 0}$ and $p_{s_i \hat{} 1}$
be extensions of $p_{s_i}$ such that $p_{s_i \hat{} 0}, p_{s_i \hat{} 1} \in D_{n+1}$,
$p_{s_i \hat{} 0} \Vdash "\underset{\sim}{y} \restriction j(n)=x_i \restriction j(n)"$
and $p_{s_i \hat{} 1} \Vdash "\underset{\sim}{y} \restriction j(n)=y_i \restriction j(n)"$.
We now define $(U_t : t\in 2^{n+1})$ as in the proof of Theorem B.
As in the proof of Theorem B, it's now easy to see that $\{ f_{k_n} : n<\omega \}$
are as required. $\square$

\textbf{\large Transferring to $L(\mathbb R)[\mathcal U]$}{\large \par}

In this section we shall prove that, assuming large cardinals, Theorem
E from the previous section also holds in $L(\mathbb R)[\mathcal U]$
where $\mathcal U$ is a selective ultrafilter on $\omega$. Suppose
that $\kappa$ is supercompact, then by the existence of an elementary
embedding $j: L(\mathbb R) \rightarrow L(\mathbb R)^{Col(\omega,\kappa)}$
(see {[}SW{]}), it follows that the results from the previous section
hold in $L(\mathbb R)$. By a result of the second author, assuming
the existence of a supercompact cardinal, any selective ultrafilter
$\mathcal U$ on $\omega$ is $\mathcal{P}(\omega)/fin$-generic over
$L(\mathbb R)$ (see {[}FA{]}). We shall now prove the main result
of this section:

\textbf{Theorem G: }Suppose that I or II hold where:

I. There exists a supercompact\textbf{ }cardinal and $\mathcal U$
be a selective ultrafilter on $\omega$.

II. $L(\mathbb R)$ is Solovay's model and $\mathcal U$ is $\mathcal{P}(\omega)/fin$-generic
over $L(\mathbb R)$.

Then in $L(\mathbb R)[\mathcal U]$, if $(g_n : n<\omega)$ is a sequence
of continuous functions from a set $\mathcal A \subseteq [\omega]^{\omega}$
to $\{0, 1\}$ and $\mathcal H \in L(\mathbb R)[\mathcal U]$ is a
selective coideal on $\omega$, then one of the following holds:

a. For every $M_0 \in \mathcal H$, there is an infinite $X\subseteq M_0$
such that the closure of $\{ g_n : n<\omega\}$ has cardinality $>2^{\aleph_0}$.

b. There is some $M\in \mathcal H$ such that $(g_n : n\in M)$ is
pointwise convergent.

\textbf{Proof: }We shall first prove the theorem under the assumption
that $L(\mathbb R)$ is Solovay's model. The proof for the supercompact
case will follow by the existence of an elementary embedding of $L(\mathbb R)$
into Solovay's model.

Find a $G_{\delta}$ set $\mathcal A^*$ such that $\mathcal A \subseteq \mathcal A^*$
and each $g_n$ extends to a continuous function $f_n$ on $\mathcal A^*$.
Each of the functions $f_n$ is in $L(\mathbb R)$, and if $(f_n : n<\omega)$
has a convergent subsequence, then so does $(g_n : n<\omega)$.

Letting $\underset{\sim}{\mathcal A}$ and $\underset{\sim}{\mathcal H}$
be a $\mathcal{P}(\omega)/fin$-names for $\mathcal A$ and $\mathcal H$,
respectively, for every $N\in [\omega]^{\omega}$, define $\mathcal A_N:=\{x \in [\omega]^{\omega} : N \Vdash x\in \underset{\sim}{\mathcal A}\}$
and $\mathcal{H}_N=\{ x\in [\omega]^{\omega} : N \Vdash x\in \underset{\sim}{\mathcal H}\}$.
Given $M_0$ and $N\in \mathcal{P}(\omega)/fin$ and $M_0 \in \mathcal H$,
we shall define the derivation process above $N$ starting from $M_0$
as follows (this will be a variant of the derivation process from
the proof of Theorem B):

We let $(X_n : n<\omega)$ be as in the proof of Theorem B. We let
$N_0=N$ and $Y_0=\mathcal A^*$. For $\alpha<\omega_1$, we shal
try to choose $(N_{\alpha}, M_{\alpha}, Y_{\alpha})$ as follows:

a. $\alpha=\beta+1$: If there is some $N$ above $N_{\beta}$ and
some $M \subseteq M_{\beta}$ such that $M \in \mathcal{H}_{N}$,
and there is some $l<\omega$ such that $X_l \cap Y_{\beta} \neq \emptyset$
and $f_m$ is constant on $X_l$ for every $m\in M$, we let $N_{\alpha}=N$,
$M_{\alpha}=M$ and $Y_{\alpha}=Y_{\beta} \setminus X_l$ where $l$
is the minimal natural number with this property.

b. $\alpha$ is a limit ordinal: In this case, choose some $N_{\alpha} \in \mathcal{P}(\omega)/fin$
above all of the conditions $\{N_{\beta} : \beta<\alpha\}$ and a
pseudeo intersection $M_{\alpha}$ of $\{M_{\beta} : \beta<\alpha\}$
such that $M_{\alpha} \in \mathcal{H}_{N_{\alpha}}$.

As there are only countably many $X_n$s, there will be a first $\alpha<\omega_1$
for which we can't carry the induction, and $\alpha$ will necessarily
be a successor ordinal. Let $(M_{N}^*, Y_{N}^*)=(M_{N_{\alpha}}, Y_{N_{\alpha}})$.
We shall also denote $\mathcal{A}_{N_{\alpha}}$ by $\mathcal{A}_N^*$
and $N_{\alpha}$ by $N^*$.

We shall now consider the following two possible cases:

Case I: There is some $N\in [\omega]^{\omega}$ such that $|Y_N^* \cap \mathcal{A}_N^*|>\aleph_0$.
In this case, by the proof of Theorem B for Solovay's model, there
is a perfect set $P=\{A_x : x\in 2^{\omega}\} \subseteq \mathcal A_N^* \cap Y_N^*$
and some infinite $X=\{k_n : n<\omega\} \subseteq M_0$ such that
$f_{k_n}(A_x)=1$ iff $x(n)=1$. Now given some $N$ with this property
such that $N^* \in \mathcal U$, as $\mathcal A_N^* \subseteq \mathcal A$,
we are done.

Case II: There is some $N_0 \in [\omega]^{\omega}$ such that for
every $N_0 \leq N$, $|Y_N^* \cap \mathcal{A}_N^*| \leq \aleph_0$.
In this case, $N_0^*$ forces that $|Y_{N_0^*} \cap \underset{\sim}{\mathcal A}| \leq \aleph_0$:
Note that $L(\mathbb R) \subseteq L(\mathbb R)[\mathcal U] \subseteq L(\mathbb{R}^*)$
where $L(\mathbb{R}^*)$ is obtained by forcing with the Mathias forcing
$\mathbb{M}_{\mathcal U}$ over $L(\mathbb R)[\mathcal U]$. By {[}DiTo{]},
there is an elementary embedding $j: L(\mathbb R) \rightarrow L(\mathbb{R}^*)$
such that $j$ fixes the reals and the ordinals. In $L(\mathbb{R}^*)$,
let $M \in \mathcal{P}(\omega)/fin$ be the Mathias generic real,
then as $M$ forces the desired statement, $\mathcal A \subseteq \mathcal{A}_M$,
$\mathbb{M}_{\mathcal U}$ is ccc and $L(\mathbb R), L(\mathbb R)[\mathcal U]$
have the same reals, the result folllows.

As in the proof of Theorem B, $N_0$ forces that $(g_n : n<\omega)$
has a converging subsequence $(g_n : n\in M)$ for some $M\in \mathcal H$.

Now consider the set of all conditions that either satisfy case II
or are of the form $N^*$ for some $N$ as in case I. As this set
is dense in $\mathcal{P}(\omega)/fin$, we're done in the case of
$L(\mathbb R)$ being Solovay's model. In the case that there exists
a supercompact cardinal $\kappa$ and $\mathcal U$ is selective,
by the existence of an elementary embedding $L(\mathbb R) \rightarrow L(\mathbb R)^{Col( \omega,\ \kappa)}$
and the last argument, there is a dense subset of $\mathcal{P}(\omega)/fin$
in $L(\mathbb R)$ consisting of conditions that force one of the
two statements forced by $N^*$ and $N_0$ in cases I and II above,
as required. $\square$

\textbf{Corollary H: }Suppose there is a supercompact cardinal and
let $\mathcal U$ be a selective ultrafilter on $\omega$. In $L(\mathbb R)[\mathcal U]$,
let $(f_n : n<\omega)$ a sequence of continuous functions from $\mathcal A$
to $2^{\omega}$, then one of the following holds:

a. There is an infinite $X\subseteq \omega$ such that the closure
of $\{f_n : n<\omega\}$ has carrdinlity $>2^{\aleph_0}$.

b. $(f_n : n< \omega)$ has a converging subsequence.

\textbf{Proof: }As in the proof of Theorem D. $\square$

\textbf{Corollary I: }Theorem F holds in $L(\mathbb R)[\mathcal U]$.

\textbf{Proof: }Let $W$ be as in the proof of Theorem F, we have
to show that if $W$ is not countably tight, then clause (b) of Theorem
F holds. Assume that a counterexample to countable tightness is given
by $Z$ and $g$ as in the proof of Theorem F. By {[}Tod{]}, we may
assume wlog that every $f\in Z$ is a pointwise limit of a subsequence
of $\{ f_n : n<\omega\}$, and hence $Z$ can be regarded as a set
of reals. We now proceed as in the proof of Theorem F. Letting $X \in \mathcal{P}(\mathbb R) \cap L(\mathbb R)[\mathcal U]$
be the domain of the $f_n$s, we let $X^*$ be a $G_{\delta}$ set
containing $X$ and $(f_n^* : n<\omega) \in L(\mathbb R)$ be continuous
extensions of the $f_n$s to $X^*$. For each $N \in \mathcal{P}(\omega)/fin$,
we let $X_N=\{x : N\Vdash x\in X\}$ and $Z_N= \{ z : N\Vdash z\in Z\}$.
Given $N\in \mathcal{P}(
\omega)/fin$, we shall describe the derivation process above $N$ similarly to
the proof of Theorem G. Let $(X_n : n<\omega)$ be as in the proof
of Theorem B. For $n<\omega$, we shall ask whether there is some
countable $A_{\frac{1}{n}} \subseteq Z_N$ such that for every condition
$N'$ above $N$, $\{ x\in X_{N'} : $ for all $f\in A_{\frac{1}{n}}$,
$\frac{1}{n+1} \leq f(x)\} =\emptyset$. If there is such an $A_{\frac{1}{n}}$
for every $n<\omega$, then $N$ will force that $0$ is in the closure
of $\underset{n<\omega}{\cup}A_{\frac{1}{n}}$. If there is an $n<\omega$
for which we can't find such an $A_{\frac{1}{n}}$, then we shall
now construct an increasing sequence of countable sets of reals $(A_{\alpha} : 0<\alpha<\alpha_*)$
and an increasing sequence of conditions as follows (where $\alpha_*$
will eventually be a countable ordinal): We let $N_0=N$. At a limit
stage $\delta$, we let $A_{\delta}=\underset{\alpha<\delta}{\cup}A_{\alpha}$
and $N_{\delta}$ be an upper bound of $\{N_{\alpha} : \alpha<\delta\}$.
At a successor stage $\alpha+1$, we ask whether there is some $X_n$,
a condition $N'$ above $N_{\alpha}$ and some countable $A \subseteq Z_{N'}$
that contains $A_{\alpha}$ such that $N'$ forces that $\{x\in X\cap X_n : \epsilon \leq f(x)$
for all $f\in A\} =\emptyset$. If there are such $N'$ and $A$,
we let $A_{\alpha+1}=A$ and $N_{\alpha+1}=N'$. Necessarily, there
will be a minimal successor ordinal $\alpha_*$ for which we can no
longer proceed. Let $B=A_{\alpha_*-1}$ and $X(B)= \{x\in X_{N_{\alpha_*-1}} : \epsilon \leq f(x)$
for every $f\in A_{\alpha_*-1} \}$. By the assumption on $\epsilon$,
$N_{\alpha_*-1}$ and $B$, we get that $X(B)$ is necessarily uncoountable
and we can now repeat the argument from the proof of Theorem F. $\square$ 

$\\$

\textbf{\large References}{\large \par}

{[}Bo{]} J. Bourgain, Compact Sets of First Baire Class, Bull. Soc.
Math. Belg, 1977

{[}BFT{]} J. Bourgain, D. H. Fremlin and M. Talagrand, Pointwise Compact
Sets of Baire-Measurable Functions, Amer. J. Math. \textbf{100 }(1978),
845-886

{[}BST{]} Karen Bakke Haga, David Schrittesser and Asger Toernquist,
Maximal Almost Disjoint Families, Determinacy and Forcing, arXiv:1810.03016

{[}DiTo{]} Carlos Augusto Di Prisco and Stevo Todorcevic, Perfect-Set
Properties in $L(\mathbb R)$, Advances in Mathematics, Volume 139,
Issue 2, 10 November 1998, Pages 240-259

{[}Fa{]} Ilijas Farah, Semiselective coideals, Mathematika 45 (1998),
no. 1, 79--103

{[}Go{]} Gilles Godefroy, Compacts De Rosenthal, Pacific Journal of
Mathematics Vol. 91, No. 2 1980

{[}HS{]} Haim Horowitz and Saharon Shelah, Can You Take Toernquist's
Inaccessible Away?, arXiv:1605.02419

{[}Ma{]} A. R. D. Mathias, Happy Families, Ann. Math. Logic \textbf{12
}(1977), no. 1, 59-111.

{[}Ro{]} H. P. Rosenthal, Pointwise Compact Subsets of the First Baire
Class, Amer. J. Math. \textbf{99 }(1977), 362-378

{[}SW{]} Saharon Shelah and Hugh Woodin, Large cardinals imply that
every reasonably definable set of reals is Lebesgue measurable. Israel
J. Math. 70 (1990), no. 3, 381--394

{[}To{]} Asger Toernquist, Definability and Almost Disjoint Familes,
Advances in Mathematics 330, 61-73, 2018

{[}Tod{]} Stevo Todorcevic, Topics in Topology, Springer Lecture Notes
in Mathematics, Vol. 1652. Springer-Verlag, Berlin, 1997. viii+153
pp. ISBN: 3-540-62611-5

$\\$

haim@math.toronto.edu

stevo@math.utoronto.ca

Department of Mathematics

University of Toronto

Bahen Centre, 40 St. George St., Room 6290

Toronto, Ontario, Canada M5S 2E4
\end{document}